\documentstyle{article}
\begin{document}
\begin{titlepage}
\title{The Weyl character formula, the half -spin representations,
and equal rank subgroups}
\author{Benedict Gross\footnote{Mathematics Department, Harvard University, Cambridge, MA 02138}, Bertram Kostant\footnote{Mathematics Department, Massachussetts Institute of Technology, Cambridge, MA 02139}, Pierre Ramond\footnote{Physics Department, Uni
versity of Florida,Gainesville, 
Fl 32611}, and Shlomo Sternberg$^*$}
\maketitle
\centerline{ {\it Proc Natl. Acad. Sci. USA}, vol 95, pp. 8441-8442, July 1988}
\vskip .5cm
\centerline{\bf {Abstract}}
\vskip .5cm
Let B be a reductive Lie subalgebra of a semi-simple Lie algebra of the same rank both over the complex numbers. To each finite dimensional irreducible representation $V_\lambda$ of F we assign a multiplet of irreducible representations of B with m elemen
ts in each multiplet, where m is the index of the Weyl group of B in the Weyl group of F. We obtain a generalization of the Weyl character formula; our formula gives the character of $V_\lambda$ as a quotient whose numerator is an alternating sum of the c
haracters in the multiplet associated to $V_\lambda$ and whose denominator is an alternating sum of the characters of the multiplet associated to the trivial representation of F. 
\vfill
\begin{flushright}    UFIFT-HEP-98-21 \end{flushright}
\end{titlepage}

Spin(9) is the light cone little group of physical theories in 10+1 dimensions 
and can also be viewed as the little group of 
massive representations in 9+1 dimensions. Computations in N=1 supergravity 
in 10+1 dimensions led to the empirical discovery of families of triples
 of irreducible representations of Spin(9) with several remarkable 
properties: the dimension of one of the three is equal to the sum of the 
dimensions of the other two, and the second order Casimir, together with 
several of the higher order Casimirs take on the same value in all three 
representations (T. Pengpan and P. Ramond, unpublished data). The purpose of this 
note is to explain this phenomenon 
and place it in a general setting. In so doing, we obtain a formula 
relating virtual representations of Lie groups of the same rank which 
reduces to the Weyl character formula when the smaller group is a 
maximal torus.

So let $B\subset F$ be two  
Lie algebras over the complex numbers with $F$ 
semi-simple and $B$ reductive and having  the same rank.  Choose a 
common Cartan subalgebra $H$ with the roots of $B$ being a subset of
the roots of $F$. Thus the Weyl group $W(B)$ of $B$ is a subgroup of the Weyl 
group $W(F)$ of $F$. Choose the positive roots consistently. Then the 
positive Weyl chamber ${\cal W}_F$ of $F$ is contained in the positive 
Weyl chamber ${\cal W}_B$ of $B$. Let $C\subset W(F)$ be those elements 
in $W(F)$ which map ${\cal W}_F$ into ${\cal W}_B$. So the cardinality of 
$C$ is the index of $W(B)$ in $W(F)$ which we denote by ${\bf m}$ and
$${\cal W}_B=\bigcup_{w\in C}w{\cal W}_F,$$
while
$$W(F)=W(B)\cdot C.$$
(In the case of Spin(9), we take $B=B_4$ and $F=F_4$ in which case 
$C$ consists of three elements.)

We let $\rho_F$ denote one half the sum of the positive roots of $F$ 
and $\rho_B$ denote one half the sum of the positive roots of $B$. Let
$\lambda $ be a dominant weight of $F$. Let $V_\lambda$ be 
the corresponding irreducible representation of $F$ 
 considered as a $B$-module by restriction.
For each $c\in C$, let 
$$c\bullet \lambda:=c(\lambda+\rho_F)-\rho_B.$$
Then $c\bullet \lambda$ is a dominant weight for $B$.  We let 
$U_{c\bullet \lambda}$
denote the irreducible representation of $B$ with
 highest weight $c\bullet \lambda$. Finally, let $P$ denote the orthogonal 
complement of $B$ in $F$ (under the Killing form of $F$) so the 
adjoint representation of $B$ on $F$ gives an embedding of $B$ 
into the orthogonal algebra $o(P)$. Since $P$ is even dimensional, $o(P)$ 
has two half-spin representations. Let $S^\pm$ denote these half-spin 
representations considered as $B$ modules. We claim that 
\begin{equation}
\label{VtimesS}
V_\lambda\otimes S^+-V_\lambda\otimes S^-=\sum_{c\in C}
\hbox{sgn}(c)U_{c\bullet \lambda}
\end{equation}
(for the appropriate choice of $\pm$). Equation (\ref{VtimesS}) is to 
be regarded as an equality in the ring of virtual representations of 
$B$. 

Recall that  ${\bf m}$ denotes   the index of $W_B$ in $W_F$
. The elements $c\bullet \lambda$  as $c$ ranges over
$c\in C$ are distinct and hence any $\lambda\in {\cal W}_F$
defines a
${\bf m}$-multiplet $\{c\bullet \lambda\},\,c\in C$ in ${\cal
W}_B$. This generalizes the triplets which arise in the case
where
$B= B_4$ and $F=F_4$. The multiplets may be abstractly
characterized as follows: Using the Harish-Chandra isomorphism
there is a natural injective homomorphism $\eta:Z_F\to Z_B$
where $Z_F$ and $Z_B$ are respectively the centers of the
enveloping algebras of $F$ and $B$. Let $Z_B^F\subset Z_B$ be the
image of $\eta$. Let ${\cal W}_B^*$ be the set of all dominant
weights $\nu \in {\cal W}_B$ such that $\nu + \rho_B$ is a
regular integral weight for $F$. Define a equivalence relation
in ${\cal W}_B^*$ by putting $\mu\sim \nu$ if the infinitesmal
characters of $U_{\mu}$ and $U_{\nu}$ agree on $Z_B^F$. \vskip 1pc
{\bf Proposition} {\it The equivalence classes all have
cardinality $\bf m$ and each such class consists of the multiplets appearing on
the right side of (\ref{VtimesS}) for a suitable $\lambda$.} 
\vskip 1pc Note that if $B$ is simple,
as in the case $B=B_4$ then the quadratic Casimir is in $Z_B^F$
and hence takes the same values for the representations whose
highest weight are in a multiplet.

 As far as we know the formula (\ref{VtimesS}) is not in the literature
although Wilfred Schmid informs us that he was aware of the
result.  The equation (\ref{VtimesS}) may be rewritten as
in (\ref{chFrepintermsofquoB}) below where it takes the form of a branching law.
The proof of (\ref{VtimesS}) combines the use of the section $C$ with
Weyl's character formula. A
recent citation for the idea of using a section of Weyl groups in
combination with the character formula to obtain a branching law
appears in 
Section 8.3.4 in ref\cite{GW}. There the branching is for the pair $C_n,C_{n+1}$.

Because the two representations occurring on the left hand side of 
(\ref{VtimesS}) have the
same dimension, we conclude that 

\begin{equation}
\label{generaldimform}
\sum_{c\in C}\hbox{sgn}(c)\hbox{dim}
(U_{c\bullet \lambda})=0.
\end{equation}

To prove (\ref{VtimesS}) we examine the Weyl character formula for 
$F$ which says that the character of $V_\lambda$ is given by
$$\hbox{ch}(V_\lambda) = \frac{A^F_{\lambda+\rho_F}}{A^F_\rho}$$
where
$$A^F_{\lambda+\rho_F}=\sum_{w\in W(F)}\hbox{sgn}(w)e^{w(\lambda+\rho_F)}$$
and where the Weyl denominator has the two expressions
$$A^F_{\rho_F}=\sum_{w\in W(F)}\hbox{sgn}(w)e^{w(\rho_F)}=
\prod\left(e^{\frac{\phi}{2}}-e^{-\frac{\phi}{2}}\right)$$
where the product is over all positive roots of $F$.

We do the summation that occurs in $A^F_{\lambda+\rho_F}$ as a double sum (over 
$W(F)=W(B)\cdot C$) to get
\begin{equation}
\label{AFis=}
A^F_{\lambda+\rho_F}=
\sum_{c\in C}\hbox{sgn}(c)A^B_{c\bullet \lambda+\rho_B}
\end{equation}
where $A^B_{c\bullet \lambda+\rho_B}$ is the numerator of the 
Weyl character formula of $B$ associated to $U_{c\bullet \lambda}$. We may 
 apply this to the Weyl denominator as well, taking $\lambda=0$. This gives
\begin{equation}
\label{chFrepintermsofB}
\hbox{ch}(V_\lambda)=\frac{\sum_{c\in C}\hbox{sgn}(c)A^B_{c\bullet \lambda+\rho_B}}
{\sum_{c\in C}\hbox{sgn}(c)A^B_{c\bullet 0+\rho_B}}.
\end{equation}

We may also rewrite the product expression for the 
Weyl denominator for $F$ as follows: Write the product over all 
positive roots of $F$ as the product over all positive roots of $B$ times 
the product over all positive missing roots:
$$A^F_{\rho_F}=A^B_{\rho_B}D$$
with
$$\Delta:=\prod_{\Phi^+(F/B)}\left(e^{\frac{\psi}{2}}-e^{-\frac{\psi}{2}}\right)$$
where $\Phi^+(F/B)$ denotes the set of
 positive roots of $F$ that are not roots of $B$. 
We thus obtain
\begin{equation}
\label{generalformula}
\hbox{ch}(V_\lambda)=\frac1{\Delta}\sum_{c\in C}\hbox{sgn}(c)\hbox{ch}
(U_{c\bullet \lambda}).
\end{equation} 
Multiplying out the product in $\Delta$ we see that 
$$\Delta=\hbox{ch }S^+-\hbox{ch }S^-$$
and multiplying (\ref{generalformula}) by $\Delta$ proves (\ref{VtimesS}).

If we divide the numerator and denominator of (\ref{chFrepintermsofB}) 
by $A^B_{\rho_B}$, the Weyl denominator for $B$, we obtain
\begin{equation}
\label{chFrepintermsofquoB}
\hbox{ch}(V_\lambda)=
\frac{\sum_{c\in C}\hbox{sgn}(c)\hbox{ch}U_{c\bullet\lambda}}
{\sum_{c\in C}\hbox{sgn}(c)\hbox{ch}U_{c\bullet 0}} 
\end{equation}
and comparing (\ref{chFrepintermsofquoB}) with 
(\ref{generalformula}), or simply taking $\lambda=0$ in 
(\ref{VtimesS}),  we see that 
\begin{equation}
\label{Dassum}
\Delta=\sum_{c\in C}\hbox{sgn}(c)\hbox{ch}U_{c\bullet 0}.
\end{equation} 

Notice that $\Delta$ vanishes on the hyperplanes determined by the $\psi\in
\Phi^+(F/B)$ 
 rather than on all the root hyperplanes as in the Weyl character 
formula. So the right hand side of (\ref{generalformula}) or equivalently 
of (\ref{chFrepintermsofquoB}) makes sense on the complement of the hyperplanes
 corresponding to the $\psi\in
\Phi^+(F/B)$ .

If we take the subgroup $B$ to be the maximal torus itself, 
equations (\ref{VtimesS}), (\ref{chFrepintermsofB}), (\ref{generalformula}), 
or (\ref{chFrepintermsofquoB})
 are just restatements of the Weyl character formula. If we take $B=B_4$
and $F=F_4$, there are exactly three terms on the right of (\ref{VtimesS}), 
since the Weyl group of $B_4$ has index three in the Weyl group of $F_4$. 
Thus, in this case, equation (\ref{generaldimform}) says that the dimension
 of the representation occurring with a minus sign on the right of 
(\ref{generaldimform}) is equal to the sum of the dimensions of
 the other two. In the case $B=D_8$ and $F=E_8$ there are 135
 terms on the right of (\ref{VtimesS}).

If we take $B=D_n$ and $F=B_n$ we obtain a formula for the character 
of a representation of $o(2n+1)$ in terms of characters of $o(2n)$; 
as the index of the Weyl groups is two, the right hand side of 
(\ref{chFrepintermsofB}) contains two terms in the numerator and in 
the denominator. Explicitly, if we make the usual choice of Cartan 
subalgebra and positive roots,
$\epsilon_i\pm \epsilon_j\ (1\leq i<j\leq n)$ for 
$D_n$ and these together with $\epsilon_i, \ (1\leq i\leq n)$ for $B_n$, 
the interior of positive Weyl chamber ${\cal W}_F$ for $B_n$ consists of all 
$x=(x_1,\dots,x_n)=x_1\epsilon_1+\cdots +x_n\epsilon_n$ with 
$$x_1>x_2>\cdots >x_n>0$$
while the interior of positive Weyl chamber ${\cal W}_B$ of $D_n$ consists of all $x$ 
satisfying
$$x_1>x_2>\cdots>x_{n-1}>|x_n|.$$
Thus $C$ consists of the identity $e$ and the reflection $s$, which 
changes the sign of the last coordinate. Here
$$\rho_F=\left(\frac{2n-1}{2},\frac{2n-3}{2},\cdots, \frac12\right)$$
whereas
$$ \rho_B= (n-1,n-2,\dots,1,0).$$
Thus for $\lambda=(\lambda_1,\dots,\lambda_n)$ we have
$$e\bullet\lambda=(\lambda_1+\frac12,\dots,\lambda_n+\frac12)\ \ \
\hbox{and }\ 
s\bullet \lambda=(\lambda_1+\frac12,\dots,\lambda_{n-1}+\frac12, -\lambda_n
-\frac12).$$
Equation (\ref{chFrepintermsofB}) becomes
\begin{equation}
\label{ortho}
\hbox{ch}V_\lambda
=\frac{\hbox{ch}U_{e\bullet \lambda}-\hbox{ch}U_{s\bullet\lambda}}
{U_{e\bullet 0}-U_{s\bullet 0}}.
\end{equation}
Notice that in this special case $S^+=U_{e\bullet 0}$ 
and $S^-=U_{s\bullet 0}$ and are the actual half spin representations of 
$o(2n)$. (For $n=1$ (\ref{ortho}) reduces to the formula for the summation of a geometrical sum.)

\end{document}